\documentclass[11pt]{amsart}
\usepackage{amscd}
\usepackage{amssymb}

\allowdisplaybreaks[1]

\newtheorem{thm}{Theorem}[section]
\newtheorem{rema}[thm]{Remark}

\newtheorem{lem}[thm]{Lemma}
\newtheorem{prop}[thm]{Proposition}

\theoremstyle{definition}

\numberwithin{equation}{section}

\begin{document}


\title{The spaces $\mathrm{H}^n(\mathfrak{osp}(1|2),M)$ for some weight modules $M$}

\begin{abstract}
We entirely compute the cohomology for a natural and large class of $\mathfrak{osp}(1|2)$ modules $M$. We study the restriction to the $\mathfrak{sl}(2)$ cohomology of $M$ and apply our results to the module $M={\mathfrak D}_{\lambda,\mu}$ of differential operators on the super circle, acting on densities. 
\end{abstract}

\author{Didier Arnal, Mabrouk Ben Ammar and Bechir Dali}

\address{Didier Arnal\\
Institut de Math\'ematiques de Bourgogne\\
UMR CNRS 5584, Universit\'e de Bourgogne\\
U.F.R. Sciences et Techniques B.P. 47870\\
F-21078 Dijon Cedex, France} \email{Didier.Arnal@u-bourgogne.fr}

\address{Mabrouk Ben Ammar\\
D\'epartement de Mathematiques, Facult\'e des Sciences de Sfax,\\
Route de Soukra, km 3,5, B.P. 1171, 3000 Sfax, Tunisie.} \email{mabrouk.benammar@fss.rnu.tn}

\address{Bechir Dali\\
D\'epartement de Math\'ematiques, Facult\'e des Sciences de Bizerte,\\
7021 Zarzouna, Bizerte, Tunisie.}\email{bechir.dali@fss.rnu.tn}

\thanks{This work was supported by the CNRS-DGRSRT project 06/S 1502. D. Arnal thanks the facult\'es des Sciences de Sfax and Bizerte for their kind hospitality, M. Ben Ammar and B. Dali thank the universit\'e de Bourgogne for its hospitality.}



\date{01/07/09}

\maketitle


\section{Introduction}

The simplest Lie  superalgebra is the algebra $\mathfrak{osp}(1|2)$.
For such an algebra, the notion of Cartan subalgebra and weight module is well known
(see section 2 for definitions and notations). In this paper, we
consider such a weight module $M$, with moreover the assumption
that one of the odd element (noted here $A$) acts through a surjective map.

This generalizes the notion of $\ell\downarrow$ modules for
$\mathfrak{sl}(2)$ \cite{m}, a class of modules admitting a finite
dimensional and nontrivial extension, but our main motivation is
the study of deformations of some actions of vector fields on the
supercircle or the superspace ${\mathbb R}^{1|1}$, this theory was developped by Ovsienko and many other
authors and some conjectures about the cohomology of natural modules
coming from the action of $\mathfrak{osp}(1|2)$ on differential
operators on densities were presented (see \cite{c,bb,g}). The first cohomology group
for this module was computed by Basdouri and Ben Ammar \cite{bab}, it was
conjectured that the second cohomology group would be generated by
cup-product of nontrivial 1 cocycles, that the 2 cocycles
whose $\mathfrak{sl}(2)$ restriction is trivial are trivial, and so one.

In this paper, we first entirely determine the cohomology for our
$\mathfrak{osp}(1|2)$ module $M$ and prove that the restriction map is one to one from ${\rm H}^n(\mathfrak{osp}(1|2),M)$ to ${\rm H}^n(\mathfrak{sl}(2),M)$. Then we apply this to the
module of the differential operators on densities, computing
completely their cohomologies and explicitely describing the cocycles.


\section{Definitions and notations}

First, we define the Lie superalgebra $\mathfrak{osp}(1|2)$ and the module $M$.
We define the superalgebra $\mathfrak g=\mathfrak{osp}(1|2)$ as the real algebra
 whose basis is $(H,X,Y,A,B)$. The elements $H$, $X$ and $Y$ are even (with parity
 0, or in $\mathfrak g_0$) and the elements $A$, $B$ are odd (with parity 1,
 or in $\mathfrak g_1$), the bracket is graded antisymmetric, we denote this
 property by $$[U,V]=-(-1)^{UV}[V,U].$$ The commutation relations are :
$$\begin{array}{llll}
&[H,X]=X,~~&[H,Y]=-Y,~~&[X,Y]=2H,\\
&[H,A]=\frac{1}{2}A,~~&[X,A]=0,~~&[Y,A]=-B,\\
&[H,B]=-\frac{1}{2}B,~~&[X,B]=A,~~&[Y,B]=0,\\
&[A,A]=2X,~~&[A,B]=2H,~~&[B,B]=-2Y.
\end{array}
$$
The bracket satisfies the graded Jacobi identity
$$
(-1)^{UW}[[U,V],W]+(-1)^{VU}[[V,W],U]+(-1)^{WU}[[W,V],U]=0.
$$
We consider the subalgebra $\mathbb R H$ as the Cartan subalgebra of
$\mathfrak{osp}(1|2)$, its adjoint action is trivially split, with roots
$0$, $\pm \frac{1}{2}$, $\pm 1$.

The even subalgebra $\mathfrak g_0$ of $\mathfrak{osp}(1|2)$ is of course
the simple Lie algebra $\mathfrak{sl}(2)$. From the relations, it is clear
that, as a graded Lie algebra, $\mathfrak{osp}(1|2)$ is generated by its odd
part $\mathfrak g_1=\mathrm{Span}(A,B)$.

We consider here a special class of $\mathfrak{osp}(1|2)$ modules $M$.
We first suppose $M$ is a complex $\mathbb Z_2$ graded vector space
$M_0\oplus M_1$ (the elements of $M_i$ are said homogenous with parity $i$) and the $H$ action is diagonalized on $M$, that is we
decompose $M$ (and thus $M_0$ and $M_1$) into weight spaces $M^\alpha$
(resp. $M_i^\alpha$) :
$$
M=\bigoplus_{\alpha\in\Sigma}M^\alpha,\quad Hv_\alpha=\alpha v_\alpha,~~\forall v_\alpha\in M^\alpha.
$$
($\Sigma\subset\mathbb C$ is the set of weights).

If $V$ is a $H$-invariant vector subspace, then $V$ itself can be decomposed
in $V=\bigoplus_{\alpha \in \Sigma}V^\alpha$ with $V^\alpha=M^\alpha\cap V$. For instance, each $M_i$ can be decomposed.

The commutation relations imply directly
$$\begin{array}{lllll}
A M_i^\alpha&\subset M_{i+1}^{\alpha+\frac{1}{2}},~~X M_i^\alpha\subset M_i^{\alpha+1},\\
B M_i^\alpha&\subset M_{i+1}^{\alpha-\frac{1}{2}},~~Y M_i^\alpha\subset M_i^{\alpha-1}.
\end{array}
$$
Then we add the condition that the action of $A$ is onto (or equivalently $X$ is onto).
This conditions implies that $M$ does not have any minimal weight vector
$v$, with weight $\alpha_0$. Indeed, if such a vector exists, the relation
$v=Aw=A\sum_{\beta\in\Sigma}w_\beta$ ($w_\beta \in M^\beta$) implies
$$\begin{array}{l}
Hv=\alpha_0 v=\sum_\beta HAw_\beta=\sum_\beta(\beta+\frac{1}{2})
Aw_\beta=\sum_\beta \alpha_0Aw_\beta,
\end{array}$$
or $Aw_\beta=0$ if $\beta\neq\alpha_0-\frac{1}{2}$, and $0\neq v=
Aw_{\alpha_0-\frac{1}{2}}$, therefore $\alpha_0-\frac{1}{2}\in\Sigma$, which is impossible. Then our modules
 $M$ are infinite dimensional.

For $\mathfrak{sl}(2)$, the simple modules for which $X$ are onto are the modules $\ell\downarrow$. It is well known that these modules are the only (with the `symmetric' case $\ell\uparrow$) $\mathfrak{sl}(2)$-modules admitting finite dimensional nontrivial extensions for some values of $\ell$ (see \cite{m}).

We now consider the cohomology groups $\mathrm{H}^n(\mathfrak{osp}(1|2),M)$ of these modules.
A $n$ cochain is a mapping $f$ from $\mathfrak{osp}(1|2)^n$ to $M$ which is $n$ linear and graded antisymmetric:
$$
f(U_1,\dots,U_i,\dots,U_j,\dots,U_n)=-(-1)^{U_iU_j}f(U_1,\dots,U_j,\dots,U_i,\dots,U_n).
$$
Defining the graded sign $\varepsilon_U(\sigma)$ for a permutation $\sigma\in\mathfrak S_n$ acting on the elements $U_i$ as the product $\varepsilon(\sigma)\varepsilon(\tau)$ of the usual sign $\varepsilon(\sigma)$ of $\sigma$ by the sign of the induced permutation $\tau$ on the set of indices $i$ for odd elements $U_i$, we have :
$$
f(U_{\sigma(1)},\dots,U_{\sigma(n)})=\varepsilon_U(\sigma)f(U_1,\dots,U_n).
$$
Due to this property, we use the following notation:
$$
U_1\cdots U_n=\frac{1}{n!}\sum_{\sigma\in\mathfrak S_n}\varepsilon_U(\sigma)(U_{\sigma^{-1}(1)}\otimes\dots\otimes U_{\sigma^{-1}(n)})
$$
and for any $\sigma\in\mathfrak S_n$, 
$$
f(U_1\dots U_n)=\varepsilon_U(\sigma)f(U_{\sigma^{-1}(1)}\otimes\dots\otimes U_{\sigma^{-1}(n)})=\varepsilon_U(\sigma)f(U_{\sigma^{-1}(1)},\dots, U_{\sigma^{-1}(n)}).
$$

The cochain $f$ is homogeneous with parity $f$ if $f(\mathfrak g_{i_1}\otimes\dots\otimes \mathfrak g_{i_n})\subset M_{f+\sum i_j}$. The space of $n$ cochain is denoted $C^n(\mathfrak{osp}(1|2),M)$, or $C^n$ if no confusion is possibe.

On such a cochain $f$, the coboundary operator is defined, using the Koszul rule for signs, by the relation (see \cite{k} for instance):
$$\aligned
&(\partial f)(U_0,\dots,U_n)=\sum_{i=0}^n (-1)^i(-1)^{U_i(f+U_0+\cdots +U_{i-1})}U_i f(U_0,\dots,\hat{\imath},\dots,U_n)+\\
&~~+ \sum_{0\leq i<j\leq n}(-1)^{i+j}(-1)^{U_i(U_0+\dots+U_{i-1})}(-1)^{U_j(U_0+\cdots+\hat{\imath}+\cdots+U_{j-1})}f([U_i,U_j],U_0,\dots,\hat{\imath},\dots,\hat{\jmath},\dots, U_n).
\endaligned
$$

If $f$ is a $n$ cochain, $\partial f$ is a $n+1$ cochain with the same parity $f$, we can verify directly that $\partial\circ \partial=0$ (or we can use a shift on degree and usual cohomology computations). The $n$ cocycles are the $n$ cochains such that $\partial f=0$, the $n$ coboundaries are the cochains in the image of $\partial$, we put as usual
$$
Z^n=\ker(\partial : C^n\longrightarrow C^{n+1}),~~~~B^n=\partial(C^{n-1}),~~~~\mathrm{H}^n(\mathfrak{osp}(1|2),M)=Z^n/B^n.
$$
$\mathrm{H}^n(\mathfrak{osp}(1|2),M)$ is the $n^{\text{th}}$ cohomology group for the module $M$.


\section{The cohomology}


\subsection{$\mathfrak{osp}(1|2)$ cohomology}


The cohomology is described by the following

\begin{thm}\label{coho}{\rm (The groups $\mathrm{H}^n(\mathfrak{osp}(1|2),M)$)}

Let us denote by $\ker A$ (respectively $\ker B$) the subspaces of $M$, kernel of the morphism $v\mapsto Av$ (respectively $v\mapsto Bv$) in $M$. Then we have the following linear isomorphisms:
\begin{itemize}
\item[\sl(i)]  $\mathrm{H}^0(\mathfrak{osp}(1|2),M)=\ker A\cap \ker B$.
\item[\sl(ii)]  $\mathrm{H}^1(\mathfrak{osp}(1|2),M)\simeq(\ker A\cap\ker B)\oplus ((\ker A)^{-\frac{1}{2}} /B((\ker A)^0))$.
\item[\sl(iii)]  $\mathrm{H}^2(\mathfrak{osp}(1|2),M)\simeq(\ker A)^{-\frac{1}{2}}/B((\ker A)^0)$.
\item[\sl(iv)] $\mathrm{H}^n(\mathfrak{osp}(1|2),M)=0$ if $n>2$.
\end{itemize}

\end{thm}
The realization of these isomorphisms will be explicitly  detailed  in the proof.
Before to prove this theorem, we shall give some preliminary results.

First we say that a $n$ cochain $f$ is \underline{reduced} if $f(AU_2\cdots U_n)=0$ for any $U_2,\dots,U_n$ in $\{A,H,B,Y\}$. Observe that if $f$ is reduced then we have also $f(XU_2\cdots U_n)=0$ for any $U_2,\dots,U_n$ in $\{A,H,B,Y\}$, since 
$$
0=(\partial f)(A^2U_2\cdots U_n)=-f([A,A]U_2\cdots U_n)=-2f(XU_2\cdots U_n).
$$

\begin{prop} {\rm(Each cochain is cohomologous to a reduced one)}

Let $f$ be a $n$ cochain. Then there exists a $n-1$ cochain $g$ such that $f-\partial g$ is reduced.

\end{prop}
\noindent
{\bf Proof.}
If $n=0$, any cochain is reduced and there is nothing to do. Suppose now $n>0$, Recall that the vectors $X,A,H,B,Y$ are root vectors with respective weight $1,\frac{1}{2},0,-\frac{1}{2},-1$. Define the weight of $U_1\otimes\dots \otimes U_n$ as the sum of the weights of the vectors $U_i$.

First, we kill  $f(A^n)$. Indeed, if $g_0$ is the $n-1$
cochain such that $g_0(U_1\dots U_{n-1})=0$ except if
$U_1=\dots=U_{n-1}=A$ and $g_0(A^{n-1}) =v$ where  $v$ is such that
$nAv=(-1)^f f(A^n)$
then $(\partial g_0)(A^n)=f(A^n)$ and $f_0=f-\partial g_0$ vanishes on
 $A^n$. If $n=1$, the proposition is proved.

Now, by induction, we suppose there is $g_k$ such that $f_k=f-\partial g_k$
 vanishes on any product of the form $A^{n-k}U_{n-k+1}\cdots U_n$
  with $U_j\in\{A,H,B,Y\}$.

Suppose $n-k>1$ and consider a $k+1$ product of the form $U_{n-k}\cdots U_n$.
If one of the $U_i$ is $A$, $f_k$ vanishes on $A^{n-k-1} U_{n-k}\cdots U_n$,
if there is no such $U_i$, but if $U_i=U_j=H$, then $f_k$ vanishes on
$A^{n-k-1}U_{n-k}\cdots U_n$. The monomial $T$ with maximal weight for
which $f_k(A^{n-k-1} T)$ could be not zero is thus $H B^k$ and its weight is $w(T)=-\frac{k}{2}$.

By induction, we can suppose $f'_k$ vanishes for any monomial of the form $A^{n-k}S$
and any monomial of the form $A^{n-k-1}T$ with $w(T)>\ell$ ($\ell\leq-\frac{k}{2}+\frac{1}{2}$).
We choose now $g_{k+1}^\ell(U_1\cdots U_{n-1})=0$ except if
$U_1\dots U_{n-1}=A^{n-k-2}U_{n-k}\cdots U_n$ and $w(U_{n-k}\cdots U_n)=\ell$
and $U_j\in\{H,B,Y\}$. Then for such a monomial,
$$\aligned
0&=g_{k+1}^\ell([A,A]A^{n-k-3}U_{n-k}\cdots U_n)=
g_{k+1}^\ell([A,U_j]A^{n-k-2}U_{n-k}\cdots\hat \jmath\cdots U_n)\\ &=
g_{k+1}^\ell([U_i,U_j]A^{n-k-1}U_{n-k}\cdots\hat \imath\cdots\hat \jmath\cdots U_n)\endaligned
$$
and
$$
(\partial g_{k+1}^\ell)(A^{n-k-1}U_{n-k}\cdots U_n)=
(-1)^{g_{k+1}^\ell}(n-k-1)Ag_{k+1}^\ell(A^{n-k-2}U_{n-k}\cdots U_n).
$$
We can then choose the value of $g_{k+1}^\ell$ such that
$f'_k-\partial g_{k+1}^\ell$ vanishes on any monomial of
the form $A^{n-k-1}T$ with $w(T)\geq\ell$.

By induction, we prove there is $g'$ such that $f'=f-\partial g'$ is reduced.

\begin{prop}{\rm(Localization for cocycles)}

Suppose than $f$ is a $n$ reduced cocycle. Then
\begin{itemize}
\item[\sl(i)] If $n>0$, $f=0$ if and only if $f(B^n)=0$.
\item[\sl(ii)] If $n>1$, any reduced cocycle vanishing on
$HB^{n-1}$ is a coboundary.
\end{itemize}

\end{prop}
\noindent
{\bf Proof.}
{\sl(i)} With the antisymmetry condition on $f$, the only possibly
non vanishing terms for $f$ are monomials containing $B$
(as odd vector) and $H$ and $Y$ as even vector, but each of them at most one time.
$f(U_1\cdots U_n)=0$ except if $U_1\cdots U_n$ is $HB^{n-1}$ or $B^n$ or $YB^{n-1}$ and, if $n>1$, $HYB^{n-2}$.

Now, the cocycle relation allows us to compute these vectors with the only knowledge of $f(B^n)$ :
$$\begin{array}{ll}
(\partial f)(AB^n)&=(-1)^fAf(B^n)+\sum_{i=1}^n(-1)^i(-1)^{(i-1)}f([A,B]B^{n-1})\\[4pt]
&=(-1)^fAf(B^n)-2nf(HB^{n-1})=0,
\end{array}
$$
and
$$\begin{array}{ll}
(\partial f)(B^{n+1})&=\sum_{i=0}^n(-1)^i(-1)^{f+i}Bf(B^n)+\sum_{0\leq i<j\leq n}(-1)^{i+j}(-1)^{i+j-1}f([B,B]B^{n-1})\\[4pt]
&=(n+1)\left[(-1)^fBf(B^n)+nf(YB^{n-1})\right]=0,\\[8pt]
(\partial f)(HB^n)&=Hf(B^n)+\sum_{i=1}^n(-1)^i(-1)^{f+i-1}Bf(HB^{n-1})+\\[4pt]
&\hskip 2cm+\sum_{j=1}^n(-1)^{j}(-1)^{j-1}f([H,B]B^{n-1})+\\[4pt]
&\hskip2cm+\sum_{1\leq i<j\leq n}(-1)^{i+j}(-1)^{i+j}f([B,B]HB^{n-2})\\[4pt]
&=(H+\frac{n}{2}id)f(B^n)-n\left[(-1)^fBf(HB^{n-1})+(n-1)f(HYB^{n-2})\right]=0.
\end{array}
$$
Thus a reduced cocycle $f$ is completely determined by the vector $f(B^n)$,
especially $f=0$ if and only if $f(B^n)=0$.

{\sl(ii)} Suppose now $f(HB^{n-1})=0$ and $n>1$. Then our computation proves
that $f(B^n)$ is in the kernel of $A$. We define a $n-1$ cochain $g$ by putting
$g(U_1\cdots U_{n-1})=0$ except for
$$
g(YB^{n-2})=\frac{1}{n(n-1)}f(B^n).
$$
Then
$$
(\partial g)(B^n)=n(n-1)g(YB^{n-2})=f(B^n).
$$
and
$$\begin{array}{lllll}
(\partial g)(AU_1\cdots U_{n-1})&=(-1)^gAg(U_1\cdots U_{n-1})+\\[4pt]
&+\sum_{j=1}^{n-1}(-1)^j(-1)^{U_j(g+1+U_2+\cdots+U_{j-1})}U_jg(AU_1\cdots\hat{\jmath}\cdots U_{n-1})+\\[4pt]
&+\sum_{j=1}^{n-1}(-1)^{j}(-1)^{U_j(U_1+\cdots+U_{j-1})}g([A,U_j]U_1\cdots\hat{\jmath}\cdots U_{n-1})+\\[4pt]
&+\sum_{1\leq i<j\leq n-1}(-1)^{i+j}(-1)^{U_i(1+U_1+\cdots+U_{i-1})+ U_j(1+U_1+\cdots\hat{\imath}\cdots+U_{j-1})}\\[4pt]
&\hskip 1cm g([U_i,U_j]AU_1\cdots\hat{\imath}\cdots\hat{\jmath}\cdots U_{n-1})\\[4pt]
&=\sum_{j=1}^{n-1}(-1)^{j}(-1)^{U_j(U_1+\cdots+U_{j-1})}g([A,U_j]U_1\cdots\hat{\jmath}\cdots U_{n-1}).
\end{array}
$$
This is non vanishing only if $[A,U_j]=B$ and $U_1\cdots\hat{\jmath}\cdots U_{n-1}= YB^{n-3}$ or $[A,U_j]=Y$ and $U_1\cdots\hat{\jmath}\cdots U_{n-1}=B^{n-2}$.
In the first case, we are computing $\partial g(AY^2B^{n-3})$, but, the antisymmetry condition on $\partial g$ gives $\partial g(AY^2B^{n-3})=0$. In the second case there is no such $U_j$. Thus $\partial g$ vanishes on any monomial $AU_1\dots U_{n-1}$.

Now $f-\partial g$ is a cocycle vanishing on any $AU_2\cdots U_n$ and on $B^n$, thus $f=\partial g$.\\

\noindent
{\bf Proof of Theorem \ref{coho}.}

{\sl(i)} If $n=0$, there is no coboundaries, the cocycles are the vector $f\in M$ such that $(\partial f)(U)=(-1)^{fU}Uf=0$ for any $U$ in $\mathfrak{osp}(1|2)$ these vectors are in $\ker A\cap \ker B$. Conversely, since $A$ and $B$ generate $\mathfrak{osp}(1|2)$ as an algebra, each vector in $\ker A\cap\ker B$ is 0 cocycle.

{\sl(ii)} Suppose $n=1$. We saw that up to a coboundary, $f$ is vanishing on $A$ and $X$. Thus $f(H)$ belongs to $\ker A$ since
$$
(\partial f)(AH)=(-1)^{f}Af(H)=0.
$$
Let us now decompose $f(H)$ on weight vectors :
$$\begin{array}{l}
f(H)=\sum_{\alpha \in \Sigma}v_\alpha,\qquad Hv_\alpha=\alpha v_\alpha,~~Av_\alpha=0.
\end{array}$$
Put $g=\sum_{\alpha\neq0}\frac{1}{\alpha}v_\alpha$. Then $(\partial g)(A)=(-1)^{g}Ag=0$ and $(\partial g)(H)=Hg=\sum_{\alpha\neq0}v_\alpha$. The 1 cocycle $f'=f-\partial g$ is reduced and satisfies $f'(H)\in \ker A\cap \ker H$. Now
$$\aligned
0=(\partial f')(HB)&=Hf'(B)-(-1)^fBf'(H)+\frac{1}{2}f'(B)
&=(H+\frac{1}{2}id)f'(B)-(-1)^fBf'(H).
\endaligned
$$
The first term is in $\bigoplus_{\alpha \neq-\frac{1}{2}} M^\alpha$, the second one in $B(\ker H)\subset M^{-\frac{1}{2}}$. Thus these two terms vanish. Therefore $f'(H)$ is in $\ker A\cap \ker B$. We now suppose $f(H)\in \ker A\cap\ker B$. Then $(H+\frac{1}{2}id)f(B)=0$.

On the other hand, we have $Af(B)=2(-1)^ff(H)$. Thus $f(B)$ is in the affine space of solutions for these two last equations. The corresponding linear space is $(\ker A)^{-\frac{1}{2}}$. But we can still add a coboundary $\partial g$ to $f$ with $Ag=Hg=0$, then $f(B)$ becomes $f(B)+(-1)^gBg$. That means, we can impose to look for solution in an affine space parallel to $(\ker A)^{-\frac{1}{2}}/B(\ker A\cap\ker H)$.

To be more precise, let us choose a supplementary space $V$ to $(\ker A)^{-\frac{1}{2}}$ in $M^{-\frac{1}{2}}$ and a supplementary space $W$ to $B((\ker A)^{0})$ in $(\ker A)^{-\frac{1}{2}}$:
$$
M^{-\frac{1}{2}}=(\ker A)^{-\frac{1}{2}}\oplus V=B((\ker A)^{0})\oplus W\oplus V.
$$
Up to a coboundary, $f(H)$ belongs to $\ker A\cap\ker B$ and $f(B)$ to $W\oplus V$. Write $f(B)=w+v$, we get $Af(B)=Av= 2(-1)^ff(H)$. This relation characterizes $v$ since $A|_V$ is one-to-one. We associate to $f$ the vector $(f(H),w)$ in $(\ker A\cap\ker B)\oplus W.$

Conversely, let $u$ be in $\ker A\cap \ker B$, homogeneous with parity $u$ and $v$ the unique vector in $V_{u+1}$ such that $Av=2(-1)^uu$. Choose any $w$ in $W$ and define a map $f:\mathfrak{osp}(1|2)\longrightarrow M$ by putting $f(A)=f(X)=0$, $f(H)=u$, $f(B)=v+w$, and $f(Y)=-(-1)^{f}B(v+w)$.
Then we verify directly that
$$\aligned
(\partial f)(AX)&=(\partial f)(AA)=(\partial f)(AH)=0\\
(\partial f)(AB)&=(-1)^fA(v+w)-2u=(-1)^fAv-2u=0,\\
(\partial f)(AY)&=-AB(v+w)-(v+w)=-(AB+BA)(v+w)-(v+w)\\&=-(2H+id)(v+w)=0.
\endaligned
$$
The map $\partial f$ is then a reduced 2 cocycle and moreover, we have $$\partial f(B^2)=(-1)^f2Bw+2(-(-1)^fBw)=0.$$ Thus $\partial f=0,$ that is, $f$ is a 1 cocycle.

Now, suppose $f$ is a coboundary, then there is $g$ such that
$$
Ag=0~~\text{ and }~~f(H)=u=Hg.
$$
This implies $H^2g=Hu=0$,  thus $g\in (\ker A)^0$ and $u=0$, thus $v=0$ and $f(B)=w=(-1)^gBg\in B((\ker A)^0)\cap W$, thus $w=0$. Conversely, if $v=0$ and $w=(-1)^gBg$ with $g\in(\ker A)^0$, then $f'=f-\partial g$ is a reduced 1 cocycle such that $f'(B)=0$, thus $f'=0$ and $f$ is a coboundary. Thus, the map $f\mapsto(u,w)$ realizes an isomorphism between $\mathrm{H}^1(\mathfrak{osp}(1|2),M)$ and $(\ker A\cap\ker B)\oplus W.$

We proved {\sl(ii)} since $W\simeq(\ker A)^{-\frac{1}{2}} /B((\ker A)^0).$

{\sl(iii)} Suppose $n\geq2$ and $f$ is a reduced $n$ cocyle. Since
$$
0=(\partial f)(AHB^{n-1})=(-1)^fAf(HB^{n-1}),
$$
we get as above: $f(HB^{n-1})$ is in $\ker A$.

We decompose $f(HB^{n-1})=\sum v_\alpha$ with $(H-\alpha id)v_\alpha=Av_\alpha=0$. Define the $n-1$ cochain $g$ by $g(U_1\dots U_{n-1})=0$ except for $g(B^{n-1})$ and $g(YB^{n-2})$ and
$$\begin{array}{lll}
g(B^{n-1})=\sum_{\alpha\neq-\frac{n-1}{2}}\frac{1}{\alpha+\frac{n-1}{2}}v_\alpha,\quad (-1)^gAg(YB^{n-2})=g(B^{n-1}).
\end{array}$$
Then $\partial g$ is a $n$ cocycle, the only non vanishing terms in $\partial g(AU_2\cdots U_n)$ are $Ag(YB^{n-2})$ and $g([A,Y]B^{n-2})$. Both happen only if $U_2\dots U_n=YB^{n-2}$ and
$$
(\partial g)(AYB^{n-2})=(-1)^gAg(YB^{n-2})-g(B^{n-1})=0.
$$
Thus $f'=f-\partial g$ is a reduced $n$ cocycle and $f'(HB^{n-1})=v_{-\frac{n-1}{2}}\in (\ker A)^{-\frac{n-1}{2}}$. From now on, we suppose $f$ is a reduced $n$ cocycle such that $f(HB^{n-1})\in (\ker A)^{-\frac{n-1}{2}}$.

Suppose now $n=2$.

If $f(HB)$ is in $B(\ker A\cap\ker H)$, we put $g(X)=g(A)=0$ and $(-1)^gBg(H)=f(HB)$ with $Ag(H)=Hg(H)=0$, then we choose $g(B)$ such that $Ag(B)=(-1)^g2g(H)$ and $g(Y)$ such that $Ag(Y)=(-1)^gg(B)$. Then $f-\partial g$ is a 2 cocycle vanishing on $AX$, $AA$, $AB$ and $AY$ and on $HB$. We saw that $f-\partial g$ is then a coboundary. Thus, $f$ is a coboundary.

Conversely, let $w$ be a vector in $(\ker A)^{-\frac{1}{2}}/B((\ker A)^0)$ (or in the supplementary space $W$ for $B((\ker A)^0)$ in $(\ker A)^{-\frac{1}{2}}$). Then
$$\begin{array}{llll}
ABw&=-w=(AB+BA)w,\\
-2Bw&=2HBw=(AB+BA)Bw,\\
AB^2w&=-2Bw-BABw=-Bw.
\end{array}
$$

We put $f(XU)=f(AU)=0$, for any $U$, $f(HB)=w$, put $f(B^2)=-4(-1)^fBw$, $f(HY)=-(-1)^fBw$, and $f(YB)=2B^2w$. The 3 cocycle $\partial f$ vanishes on $A^2U$ for any $U$, we consider it on $AHB$, $AB^2$, $AYB$ and $AYH$.
$$\begin{array}{llll}
(\partial f)(AHB)&=(-1)^fAw=0,\\
(\partial f)(AB^2)&=-4ABw-2f([A,B]B)=4w-4w=0,\\
(\partial f)(AHY)&=-ABw+f([A,Y]H)=w-w=0,\\
(\partial f)(AYB)&=(-1)^fAf(YB)-f([A,Y]B)+f([A,B]Y)\\
&=(-1)^f\left[2AB^2w+4Bw-2Bw\right]=0.\\
\end{array}
$$
$\partial f$ is a reduced 3 cocycle, we moreover have
$$
(\partial f)(B^3)=(-1)^f3Bf(B^2)-3f([B,B]B)=-12B^2w+6f(YB)=0.
$$
Thus $\partial f=0$, $f$ is then a reduced 2 cocycle. Now if $f=\partial g$, then 
$$
\begin{array}{l}w=f(HB)=\partial g(HB)=\left(H+\frac{1}{2}id\right)g(B)-Bg(H).
\end{array}
$$
Let $g(B)=\sum_\alpha u_\alpha$, $g(H)=\sum_\alpha x_\alpha$  and $g(A)=\sum_\alpha y_\alpha$ where $u_\alpha$, $x_\alpha$ and $y_\alpha$ are in $ M^\alpha$, then we get 
$$
\begin{array}{l}w=\sum_{\alpha\neq-\frac{1}{2}}\left((\alpha+\frac{1}{2})u_\alpha-Bx_{\alpha+\frac{1}{2}}\right)-Bx_0.\end{array}
$$
But $w$ is in $W$, thus, $(\alpha+\frac{1}{2})u_\alpha-Bx_{\alpha+\frac{1}{2}}=0$ if $\alpha\neq-\frac{1}{2}$ and then $w=-Bx_0$.
Moreover, we have
$$
\begin{array}{l}0=f(HA)=(H-\frac{1}{2}id)g(A)-Ag(H)=\sum_{\alpha\neq \frac{1}{2}}\left((\alpha-\frac{1}{2})y_\alpha-Ax_{\alpha-\frac{1}{2}}\right)-Ax_0.\end{array}
$$
Thus, $Ax_0=0$, therefore $x_0\in (\ker A)^0$ and $w=-Bx_0\in W\cap B((\ker A)^0)=\{0\}$, this implies  $f=0$.

We proved the point {\sl(iii)}.

{\sl(iv)} Suppose $n>2$. We saw that any $n$ cocycle $f$ can be choosen such that $f$ is reduced and $f(HB^{n-1})\in(\ker A)^{-\frac{n-1}{2}}$. We define $g$ by $g(U_1\dots U_{n-1})=0$ except
$$
g(HYB^{n-3})=-\frac{1}{(n-1)(n-2)}f(HB^{n-1})
$$
and $g(YB^{n-2})$, choosen such that
$$
Ag(YB^{n-2})-2(n-2)g(HYB^{n-3})=0.
$$

Then $(\partial g)(AHB^{n-2})=0$ and if $U_2\dots U_n\neq HB^{n-2}$, then the only non vanishing terms in $(\partial g)(AU_2\dots U_n)$ have the form $\pm g([A,U_j]U_2\dots\hat{\jmath}\dots U_n)$ with $[A,U_j]=Y$, which is impossible, or $[A,U_j]=B$, this means $U_j=Y$, but there is another index $i\neq j$ with $U_i=Y$ and this is still impossible or $[A,U_j]=H$, this means $U_j=B$ and $U_2\dots U_n= HB^{n-2}$, which is impossible. Thus $f-\partial g$ is reduced and vanishes on $HB^{n-1}$, it is a coboundary, $f$ is a coboundary, $\mathrm{H}^n(\mathfrak{osp}(1|2),M)=0$.

\subsection{Restriction to $\mathfrak{sl}(2)$}


We keep our notations.

\begin{lem}{\rm (Characterization for $B((\ker A)^0)$)}

Let $w$ be a vector in $M$ such that $w\in(\ker A)^{-\frac{1}{2}}$ and $Bw\in Y((\ker X)^0)$. Then $w$ is in $B((\ker A)^0)$.\\
\end{lem}

\noindent
{\bf Proof.}
We suppose $Bw=B^2v$, with $Hv=A^2v=0$. Thus $ABv+BAv=0$ and
$$
2HBv=-Bv=AB^2v+BABv,\qquad ABw=AB^2v=-Bv-BABv.
$$
But
$$
2Hw=(AB+BA)w=ABw=-w.
$$
Or $w=B(v+ABv)$. But
$$
2HAv=(AB+BA)Av=ABAv=Av.
$$
Finally:
$$
A(v+ABv)=Av+A^2Bv=Av-ABAv=0.
$$
This proves our lemma.

Let $f$ be a $n$ cochain for $\mathfrak{osp}(1|2)$. Its restriction $f|_{\mathfrak{sl}(2)}$ to $\mathfrak{sl}(2)^n$ is a $n$ cochain for the $\mathfrak{sl}(2)$ module $M$. If $f$ is a cocycle (resp. a coboundary), $f|_{\mathfrak{sl}(2)}$ is a cocycle (resp. a coboundary). The map $f\mapsto f|_{\mathfrak{sl}(2)}$ defines a map $\varphi$ from ${\rm H}^n(\mathfrak{osp}(1|2),M)$ to ${\rm H}^n(\mathfrak{sl}(2),M)$.

\begin{prop}{\rm(Restriction of $\mathfrak{osp}(1|2)$ cocycle and triviality)}

A $n$ cocycle $f$ for $\mathfrak{osp}(1|2)$ is a coboundary (a trivial cocycle) if and only if its restriction $f|_{\mathfrak{sl}(2)}$ is a $\mathfrak{sl}(2)$ coboundary. Or: $\varphi$ is one to one.
\end{prop}

\noindent
{\bf Proof.}
We just consider $n\leq2$ and $f$ choosen as in Theorem \ref{coho}.

A 0 cocycle is a vector $f$ in $\ker A\cap\ker B$, it is trivial if and only if $f=0$.

A 1 cocycle is cohomologous to a cocycle $f$ such that:
$$\aligned
f(A)&=f(X)=0,\quad
f(H)=u\in (\ker A)^0,\\
f(B)&=v+w\in V\oplus W,\quad
f(Y)=-(-1)^fB(v+w).
\endaligned
$$
Here $V$ is a supplementary space for $(\ker A)^{-\frac{1}{2}}$ in $M^{-\frac{1}{2}}$, $W$ a supplementary space for $B((\ker A)^0)$ in $(\ker A)^{-\frac{1}{2}}$, and $v$ is the only vector in $V$ such that $Av=2(-1)^fu$. We saw that $f$ is characterized by $u$ and $w$.

Suppose there is $g$ in $M$ such that $(f-\partial g)|_{\mathfrak{sl}(2)}$ vanishes, thus $Hg=f(H)=u$, since $u$ is in $M^0$, this relation forces $u=0$, therefore $v=0$. Now $Xg=f(X)=0$, $Hg=f(H)=0$ and $Yg=f(Y)=-(-1)^fBw$. Our lemma says that $w$ is in $B((\ker A)^0)$, thus $w=0$, $f=0$.\\

A 2 cocycle is cohomologous to a cocycle $f$ such that:
$$\begin{array}{lllll}
f(AU)&=f(XU)=0,\quad
f(HB)&=w\in W,\quad
f(BB)&=-4(-1)^fBw,\\[4pt]
f(HY)&=-(-1)^fBw,\quad
f(YB)&=2B^2w.\quad&~~
\end{array}
$$
And $f$ is characterized by $w$.

Suppose there is $g$ in $C^1(\mathfrak{sl}(2),M)$ such that $(f-\partial g)|_{\mathfrak{sl}(2)}$ vanishes, put:
$$\begin{array}{ll}
g(X)=\sum_\alpha x_\alpha,\quad
g(H)=\sum_\alpha h_\alpha,\quad
g(Y)=\sum_\alpha y_\alpha,\quad(x_\alpha,\,h_\alpha,\,y_\alpha\in M^\alpha).
\end{array}$$
We get
$$\begin{array}{ll}
f(XH)&=Xg(H)-Hg(X)-g([X,H])=\sum_{\alpha\neq1}((-\alpha+1)x_\alpha+Xh_{\alpha-1})+Xh_0=0,\\[4pt]
f(HY)&=Hg(Y)-Yg(H)-g([H,Y])=\sum_{\gamma\neq-1} ((\gamma+1)y_\gamma-Yh_{\gamma+1})-Yh_0=-(-1)^fBw.
\end{array}
$$
Since $Bw$ is in $M^{-1}$, this implies $Hh_0=Xh_0=0$ and $Yh_0=(-1)^fBw$. Our lemma says that $w$ is in $B((\ker A)^0)$, therefore $w=0$ and $f=0$.

\begin{rema} 
In the same way as for Theorem \ref{coho}, it is easy to compute the cohomology for the $\mathfrak{sl}(2)$ module $M$. Here it is:
$$\aligned
{\rm H}^0(\mathfrak{sl}(2),M)&=\ker X\cap \ker Y,\qquad {\rm H}^1(\mathfrak{sl}(2),M)\simeq(\ker X\cap \ker Y)\oplus (\ker X)^{-1}/Y((\ker X)^0),\\
{\rm H}^2(\mathfrak{sl}(2),M)&\simeq (\ker X)^{-1}/Y((\ker X)^0), \qquad {\rm H}^{>2}(\mathfrak{sl}(2),M)=0.
\endaligned
$$ 
\end{rema}
\section{Application to $\mathfrak D_{\lambda,\mu}$}


\subsection{Differential operators on weighted densities}


\

We define the superspace $\mathbb{R}^{1|1}$ in terms of its
superalgebra of functions, denoted by $C^\infty(\mathbb{R}^{1|1})$
and consisting of elements of the form:
\begin{equation*}
F(x,\theta)=f_0(x)+f_1(x)\theta ,
\end{equation*}
where $x$ is the even variable,  $\theta$ is  the odd variable
($\theta^2=0$) and $f_0(x),\,f_1(x)\in C^\infty(\mathbb{R})$.  We
consider the contact bracket on $C^\infty(\mathbb{R}^{1|1})$
defined on  $C^\infty(\mathbb{R}^{1|1})$ by:
\begin{equation*}\begin{array}{l}
\{F,G\}=FG'-F'G+\frac{1}{2}\eta(F)
\overline{\eta}(G),\end{array}
\end{equation*}where $\eta=\frac{\partial}{\partial
{\theta}}+\theta\frac{\partial}{\partial x}$ and
$\overline{\eta}=\frac{\partial}{\partial
{\theta}}-\theta\frac{\partial}{\partial x}$.
Let $\mathrm{Vect}(\mathbb{R}^{1|1})$ be the superspace of
 vector fields on $\mathbb{R}^{1|1}$:
\begin{equation*}\mathrm{Vect}(\mathbb{R}^{1|1})=\left\{F_0\partial_x
+  F_1\partial_\theta \mid ~F_i\in C^\infty(
\mathbb{R}^{1|1})\right\},\end{equation*} where $\partial_\theta$
stands for $\frac{\partial}{\partial\theta}$ and $\partial_x$
stands for $\frac{\partial}{\partial x} $. We can realize the algebra $\mathfrak{osp}(1|2)$ as a subalgebra of $\mathrm{Vect}(\mathbb{R}^{1|1})$:
\begin{equation*}
\mathfrak{osp}(1|2)=\text{Span}(X_1,\,X_{x},\,X_{x^2},\,X_{x\theta},\,
X_{\theta}).
\end{equation*}
where, the vector field $X_G$ is defined for any $G\in C^\infty(\mathbb{R}^{1|1})$ by
\begin{equation*}\begin{array}{l}
X_G=G\partial_x+\frac{1}{2}\eta(G)\overline{\eta}.\end{array}
\end{equation*}
Here, we have $(-X_x,X_1,-X_{x^2},2X_\theta,X_{x\theta})=(H,X,Y,A,B).$
The bracket on $\mathfrak{osp}(1|2)$ is then given by
$
[X_{F},\,X_{G}]=X_{\{F,\,G\}}.
$

We denote by $\mathfrak{F}_{\lambda}$ the space of all
 weighted densities on $\mathbb{R}^{1|1}$ of weight $\lambda$:
\begin{equation*}
\mathfrak{F}_\lambda=\left\{F(x,\theta)\alpha^\lambda~~|~~F(x,\theta)
\in C^\infty(\mathbb{R}^{1|1})\right\}\qquad(\alpha=dx+\theta d_\theta).
\end{equation*}
The action of $\mathfrak{osp}(1|2)$ on $\mathfrak{F}_{\lambda}$ is given by
\begin{equation*}
\label{superaction}
{X_G}(F\alpha^\lambda)=((G\partial_x+\frac{1}{2}\eta(G)\overline{\eta})(F)+
\lambda G'F)\alpha^\lambda.
\end{equation*}

Any differential operator $A$ on $\mathbb{R}^{1|1}$ defines a linear mapping from
$\mathfrak{F}_\lambda$ to $\mathfrak{F}_\mu$ for any $\lambda$ by: $A :F\alpha^\lambda\mapsto A(F)\alpha^\mu$,
$\mu\in\mathbb{R}$, thus, the space of differential operators
becomes a family of  $\mathfrak{osp}(1|2)$
modules denoted $\mathfrak{D}_{\lambda,\mu}$, for the natural
action:
\begin{equation*}\label{d-action}
{X_G}\cdot A={X_G}\circ A-
(-1)^{AG}A\circ {X_G}.
\end{equation*}
For more details see, for instance \cite{ab, bab, bb,g}


\subsection{Cohomology}


\

Let us consider the $\mathfrak{osp}(1|2)$-module $\mathfrak{D}_{\lambda,\mu}$ of differential operators on densities on $\mathbb{R}^{1|1}$.

We put here $p=\mu-\lambda$ and choose the following basis for $\mathfrak{D}_{\lambda,\mu}$ :
$$
a_{m,k}=x^m\partial_x^k,~~b_{m,k}=x^m\theta\partial_\theta\partial_x^k,~~c_{m,k}=x^m\theta\partial_x^k,~~d_{m,k}=x^m\partial_\theta\partial_x^k-x^m\theta\partial_x^{k+1}.
$$
(Here, $m$ and $k$ are natural integral numbers), we say that $a_{m,k}$ and $b_{m,k}$ are even vectors (see below) and $c_{m,k}$ and $d_{m,k}$ are odd vectors.

In fact they are weight vectors for the action of $H$:
$$
\begin{array}{llllll}
Ha_{m,k}&=(k-m-p)a_{m,k}\quad
&Hb_{m,k}&=(k-m-p)b_{m,k}\\
Hc_{m,k}&=(k-m-p-\frac{1}{2})c_{m,k}\quad
&Hd_{m,k}&=(k-m-p+\frac{1}{2})d_{m,k}.
\end{array}
$$

Similarly, a direct computation give the following relations for the $A$ and $B$ actions on these vectors :
$$
\begin{array}{llllll}
A~ a_{m,k}&=mc_{m-1,k},
&A~ b_{m,k}&=d_{m,k},\\
A~ c_{m,k}&=a_{m,k},
&A~ d_{m,k}&=mb_{m-1,k}
\end{array}
$$
and
$$
\begin{array}{lllll}
B a_{m,k}&=(m-2k+2p)c_{m,k}-kd_{m,k-1},
&B b_{m,k}&=d_{m+1,k}-(2\lambda+k)c_{m,k},\\
B c_{m,k}&=a_{m+1,k}+kb_{m,k-1},
&B d_{m,k}&=(m-2k+2p-1)b_{m,k}+(2\lambda +k)a_{m,k}.
\end{array}
$$
From these formulas (or directly), we can compute the $X$ and $Y$ actions, getting:
$$
\begin{array}{llll}
Xa_{m,k}&=ma_{m-1,k},&Xb_{m,k}&=mb_{m-1,k},\\
Xc_{m,k}&=mc_{m-1,k},&Xd_{m,k}&=md_{m-1,k},
\end{array}
$$
and
$$
\begin{array}{ll}
Ya_{m,k}&=(2k-2p-m)a_{m+1,k}+k(2\lambda+k-1)a_{m,k-1}+kb_{m,k-1},\\
Yb_{m,k}&=(2k-2p-m)b_{m+1,k}+k(2\lambda+k)b_{m,k-1},\\
Yc_{m,k}&=(2k-2p-m-1)c_{m+1,k}+k(2\lambda+k-1)c_{m,k-1},\\
Yd_{m,k}&=(2k-2p-m+1)d_{m+1,k}+k(2\lambda+k)d_{m,k-1}-(2\lambda+2k+1)c_{m,k}.
\end{array}
$$

From these formulas, we immediately get
$$
\ker A\cap \ker B=~\left\{\begin{array}{ll}
{\rm Span}(a_{0,0})&~\text{ if }~ p=0,\\
{\rm Span}(d_{0,k})&~\text{ if }~ p=k+\frac{1}{2},~~k\in\{0,1,2,\dots\}~\text{ and }~2\lambda+k=0,\\
0&~~\text{elsewhere.}
\end{array}\right.
$$
and
$$
(\ker A)^{-\frac{1}{2}}=~\left\{\begin{array}{ll}
{\rm Span}(a_{0,k})&~\text{ if }~ p=k+\frac{1}{2},~~k\in\{0,1,2,\dots\},\\
{\rm Span}(d_{0,k})&~\text{ if }~ p=k+1,~~k\in\{0,1,2,\dots\},\\
0&~~\text{elsewhere.}
\end{array}\right.
$$
Moreover, if $p=k+\frac{1}{2}$, $(\ker A)^0={\rm Span}(d_{0,k})$, $B((\ker A)^0)={\rm Span}(a_{0,k})$ if $2\lambda+k\neq0$, 0 if it is not the case. Similarly, if $p=k+1$, then $B((\ker A)^0)=B({\rm Span}(a_{0,k+1}))={\rm Span}(d_{0,k})$.

Now we deduce :

\begin{prop}{\rm(The cohomology for $\mathfrak{D}_{\lambda,\mu}$)}

The dimensionalities for the cohomology groups $\mathrm{H}^n(\mathfrak{osp}(1|2),\mathfrak D_{\lambda,\mu})$ are:
$$\aligned
\text{\sl(i)}&\quad \dim(\mathrm{H}^0(\mathfrak{osp}(1|2),\mathfrak D_{\lambda,\mu}))=~\left\{\begin{array}{ll}
1&~\text{ if }~\lambda=\mu,\\
1&~\text{ if }~\lambda=-\frac{k}{2}~\text{and}~\mu=\frac{k+1}{2},~~k\in\{0,1,2,\dots\},\\
0&~\text{ in the other cases}.\end{array}\right.\\
\text{\sl(ii)}&\quad\dim(\mathrm{H}^1(\mathfrak{osp}(1|2),\mathfrak D_{\lambda,\mu}))=~\left\{\begin{array}{ll}
1&~\text{ if }~ \lambda=\mu,\\
2&~\text{ if }~\lambda=-\frac{k}{2}~\text{and}~\mu=\frac{k+1}{2},~~k\in\{0,1,2,\dots\},\\
0&~\text{ in the other cases}.\end{array}\right.\\
\text{\sl(iii)}&\quad\dim(\mathrm{H}^2(\mathfrak{osp}(1|2),\mathfrak D_{\lambda,\mu}))=~\left\{\begin{array}{ll}
1&~\text{ if }~\lambda=-\frac{k}{2}~\text{and}~\mu=\frac{k+1}{2},~~k\in\{0,1,2,\dots\},\\
0&~\text{ in the other cases}.\end{array}\right.\\
\text{\sl(iv)}&\quad\dim(\mathrm{H}^n(\mathfrak{osp}(1|2),\mathfrak D_{\lambda,\mu}))=0.\endaligned
$$
\end{prop}

We refind here the results of \cite{bab} for the $\mathrm{H}^1$.

To be more precisey, in the following, we give explicit basis for these cohomology groups
\begin{itemize}
\item [\sl(i)] $\mathrm{H}^0(\mathfrak{osp}(1|2),\mathfrak D_{\lambda,\lambda})=\mathrm{Span}(id)$ and $\mathrm{H}^0\left(\mathfrak{osp}(1|2),\mathfrak D_{-\frac{k}{2}, \frac{k+1}{2}}\right)=\mathrm{Span}(\partial_\theta\partial_x^k-\theta\partial_x^{k+1})$.
\item [\sl(ii)] The space $\mathrm{H}^1(\mathfrak{osp}(1|2),\mathfrak D_{\lambda,\lambda})$ is spanned by the cohomology class of the reduced 1 cocycle $h_\lambda$ defined by:
$$
h_\lambda(X)=h_\lambda(A)=0,\quad h_\lambda(H)=-id,\quad h_\lambda(B)=\theta\cdot \quad\text{and}\quad h_\lambda(Y)=-2x \cdot.
$$
While the space $\mathrm{H}^1\left(\mathfrak{osp}(1|2),\mathfrak D_{-\frac{k}{2}, \frac{k+1}{2}}\right)$ is spanned by the cohomology classes of the reduced 1 cocycles $f_k$ and $\widetilde{f}_k$  defined respectively by:
$$
f_k(X)=f_k(A)=0,\quad f_k(H)=\partial_\theta\partial_x^k-\theta\partial_x^{k+1},\quad f_k(B)=\theta \partial_\theta\partial_x^k\,\,\text{and}\,\, f_k(Y)=2x f_k(H),
$$
$$\widetilde{f}_k(X)=\widetilde{f}_k(A)=\widetilde{f}_k(H)=0,\quad\widetilde{f}_k(B)=\partial_x^k\,\,\text{and}\,\, \widetilde{f}_k(Y)=-2k\partial_\theta\partial_x^{k-1}+2\theta(k+1)\partial_x^k.
$$
\item [\sl(iii)] A similar realization of $\mathrm{H}^2\left(\mathfrak{osp}(1|2),\mathfrak D_{-\frac{k}{2}, \frac{k+1}{2}}\right)$ is easy, we prefer to give an explicit, nontrivial, reduced 2 cocycle as a cup product. Let
$$
\Omega_k(U,V)=(f_k\vee h_{-\frac{k}{2}})(U,V):=f_k(U)\circ h_{-\frac{k}{2}}(V)-(-1)^{UV}f_k(V)\circ h_{-\frac{k}{2}}(U).
$$
Since $f_k$ and $h_{-\frac{k}{2}}$ are cocyles, a direct computation shows that $\Omega_k$ is a 2 cocycle, it is nontrivial since its restriction to $\mathfrak{sl}(2)\times\mathfrak{sl}(2)$ is nontrivial: 
$$
\Omega_k(X_f,X_g)=-(-1)^k\omega(f,g)(k\partial_\theta\partial_x^{k-1}-(k+1)\theta\partial_x^{k})
$$
where $\omega$ is the Gelfand-Fuchs cocycle defined by $\omega(f,g)=f'g''-g'f''$.
\end{itemize}

\end{document}